\documentclass[reqno,12pt]{amsart}
\usepackage[utf8]{inputenc}
\usepackage[margin=1in]{geometry}

\usepackage{amsmath, amssymb, amsfonts, latexsym, mdwlist, amsthm, amscd}
\usepackage[]{mathrsfs}
\usepackage[draft=false]{hyperref}
\hypersetup{final}
\usepackage{enumerate}

\usepackage{tikz}
\usetikzlibrary{calc,trees,positioning,arrows,chains,shapes.geometric,%
    decorations.pathreplacing,decorations.pathmorphing,shapes,%
    matrix,shapes.symbols}
\usepackage{caption}
\usepackage{subcaption}

\usepackage{graphicx}
\usepackage{tikz-cd}
\usepackage{pgfplots}
\pgfplotsset{compat=1.17}

\newtheorem{thm}{Theorem}
\newtheorem{prop}[thm]{Proposition}
\newtheorem{lem}[thm]{Lemma}

\theoremstyle{definition}
\newtheorem{defn}[thm]{Definition}

\theoremstyle{remark}
\newtheorem{rem}[thm]{Remark}

\DeclareMathOperator{\Hom}{Hom}
\DeclareMathOperator{\Ext}{Ext}
\DeclareMathOperator{\ext}{ext}
\DeclareMathOperator{\rk}{rk}
\DeclareMathOperator{\im}{im}
\DeclareMathOperator{\codim}{codim}
\DeclareMathOperator{\id}{id}

\DeclareMathOperator{\Quot}{Quot}
\DeclareMathOperator{\Fl}{Fl}
\DeclareMathOperator{\gr}{gr}

\newcommand{\ard}{A(r_1,d_1)}

\newcommand{\urdo}{\mathcal{U}^0(r_1)}
\newcommand{\oo}{{\mathcal{O}}}

\newcommand{\ignore}[1]{}

\title{Shatz strata in algebraic versal deformation spaces}
\author[]{Yinbang Lin}
\address{Tongji University, Shanghai, China} 
\email{yinbang\textunderscore lin@tongji.edu.cn}
\date{}
\subjclass[2020]{Primary: 14D20; Secondary: 14B05, 14D15}
\keywords{Algebraic versal deformation space; Shatz stratum; Harder-Narasimhan filtration; Sheaves on algebraic curves}

\begin{document}

\begin{abstract}
Over a smooth complex projective curve, we study an algebraic versal deformation space with fixed determinant of a coherent sheaf. The algebraic versal deformation space decomposes into a disjoint union of Shatz strata, namely locally closed subschemes which parametrize coherent sheaves with common Harder-Narasimhan types. We study the geometry and local topology of large unstable strata and their behavior along boundaries.
\end{abstract}
\maketitle

\section{Introduction}

Vector bundles on algebraic curves has been a central topic in algebraic geometry for decades. It is related to various areas of mathematics.

Slope stability is a basic notion for vector bundles. Imposing the stability condition, we can obtain a finite type moduli space of semistable vector bundles with fixed invariants. For an unstable vector bundle, it is well known that there is a {\em Harder-Narasimhan filtration} such that the factors are semistable. The filtration allow us to break down the vector bundle into basic building blocks.

Given a family of vector bundles parametrized by a finite type scheme $S$, we can decompose $S$ into locally closed subschemes, according to the types of Harder-Narasimhan filtrations \cite{Sha77}. We call them {\em Shatz strata}\footnote{They probably should be called Harder-Narasimhan-Shatz strata. But let's call them Shatz strata for brevity.}.
We are interested in the geometry and local topology of the Shatz strata in versal deformation spaces of coherent sheaves on complex algebraic curves. We will extend the corresponding results in \cite[\S 2]{Li97}, which focused on rank 2 cases.

Let $C$ be a smooth complex projective curve of genus $g$ and $L$ be a line bundle of degree $d$ on $C$.
We often focus on vector bundles while studying moduli problems. But we need to consider sheaves with torsion in order to allow potential applications. For example, if we consider a torsion free sheaf on a smooth surface and its restriction to a subcurve, then torsion is likely to appear.
We thus define stability (i.e. slope stability) for arbitrary coherent sheaves of positive rank.
\begin{defn}
A coherent sheaf $E$ on $C$ of rank $r>0$ is $\mu$-{\em stable} (or simply {\em stable}) if $\mu(E)<\mu(G)$ for any torsion free quotient $G$ of rank $<r$.
Semistability is defined in a similar way, replacing $<$ by $\leqslant$.
\end{defn}

Note that there are sheaves with torsion that are $\mu$-stable.

\begin{figure}[ht]
    \centering
    \includegraphics[scale=.25]{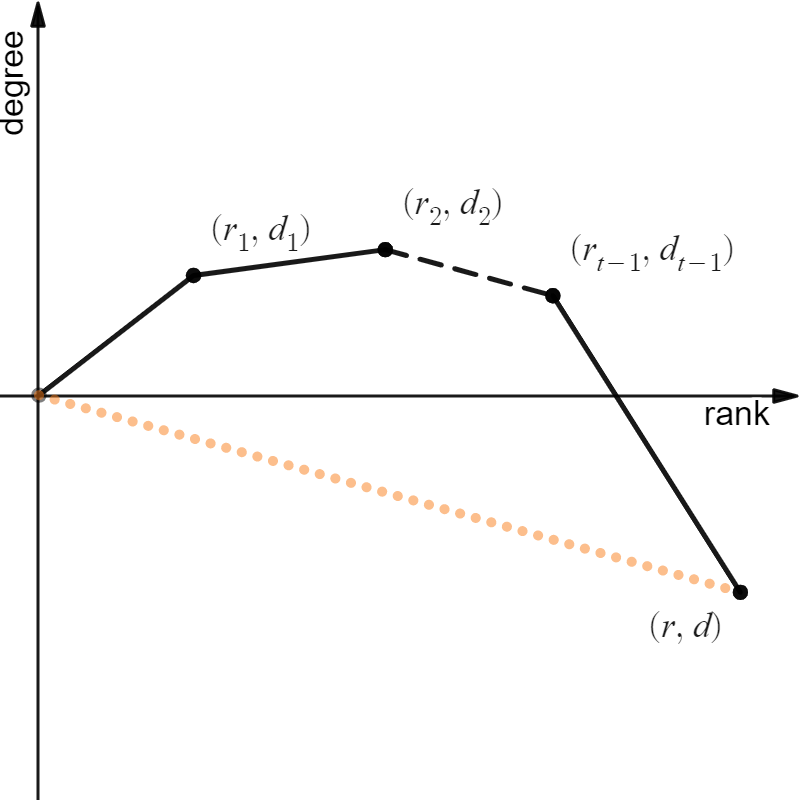}
    \caption{Harder-Narasimhan polygon}
    \label{fig:hn-polygon}
\end{figure}

For any coherent sheaf $E$ of positive rank, there is a unique {\em Harder-Narasimhan filtration} (or HN-filtration)
\[0=F_0\subsetneqq F_1\subsetneqq \dots \subsetneqq F_t=E,\]
such that
\begin{enumerate}[(a)]
    \item $F_i/F_{i-1}$ is $\mu$-semistable for $i=1,\dots,t$,
    \item $\mu(F_1/F_{0})>\mu(F_2/F_{1})>\dots>\mu(F_t/F_{t-1})$.
\end{enumerate}
Then, torsion, if any, is entirely contained in $F_1$.
Let \[\mu_{\max}(E)=\mu(F_1/F_{0})\mbox{ and }\mu_{\min}=\mu(F_t/F_{t-1}).\]
The ranks and degrees $(r_i,d_i)=(\rk F_i,\deg F_i)$ provide a polygon associated to $E$ in the rank-degree plane as in Figure~\ref{fig:hn-polygon}.
This is called the {\em Harder-Narasimhan polygon}, or HN-polygon for short.
Notice that the sheaf $E$ is unstable if and only if the HN-polygon is strictly above the segment connecting $(0,0)$ and $(r,d)=(\rk E,\deg E)$.
We say $E$ has {\em HN-type} $(r_1,d_1;r_2,d_2;\ldots;r_t,d_t)$.

\begin{defn}
Let $Z$ be a projective scheme and let $E$ be a coherent sheaf over $Z$ with $\det E \cong L$. An {\em algebraic versal deformation space} of $E$ with fixed determinant is a collection \[(B, \mathscr{E}_B;0, E),\] where $B$ is a quasi-projective scheme, $0\in B$, and $\mathscr{E}_B$ is a family of sheaves on $B\times Z$ flat over $B$ with $\det \mathscr{E}_B\cong p_Z^*L$, $p_Z\colon B\times Z\to Z$, such that
\begin{enumerate}[(a)]
    \item the restriction of $\mathscr{E}_B$ to $\{0\}\times Z$, say $\mathscr{E}_0$, is isomorphic to $E$, and that the Kodaira-Spencer map $T_0 B\to \Ext^1(\mathscr{E}_0,\mathscr{E}_0)_0$ induced by the family $\mathscr{E}_B$ is an isomorphism;
    \item For any pair of affine varieties $S_0 \subset S$ coupled with a sheaf $\mathscr{F}_S$ on $S \times Z$ flat over $S$ such that $S_0\subset S$ is closed, $\det \mathscr{F}_S\cong p_Z^*L$ and $\mathscr{F}_S|_{S_0\times Z}\cong q_Z^*E$, where $q_Z\colon S_0\times Z\to Z$, there is an analytic neighborhood $U$ of $S_0\subset S$ and an analytic map $\eta\colon (U, S_0)\to (B, 0)$ so that the restriction of $\mathscr{F}_S$ to $U\times Z$ is isomorphic to $(\eta \times \id_Z)^*\mathscr{E}_B$, extending the given isomorphism $\mathscr{E}_0\cong E$ and $\mathscr{F}_S|_{S_0\times Z}\cong q_Z^*E$.
\end{enumerate}
\end{defn}
See e.g. \cite[Definition 2.3]{Li97}. Since we only study algebraic versal deformation spaces, we often omit the word ``algebraic''.

Let $E_0$ be an unstable sheaf on $C$ with rank $r$ and determinant $L$. Suppose its HN-type is
\[(s_1,e_1;s_2,e_2;\dots;s_t=r,e_t=d).\]
Let $(A,\mathscr{E}_A;0, E_0)$ be a fixed versal deformation space of $E_0$ with fixed determinant, such that $A$ is a {\em smooth} affine variety.
We refer the reader to \cite{Li97} for the construction of a versal deformation space.
Notice that \[\dim A=\ext^1(E_0,E_0)_0=-\chi(E_0,E_0)_0+\hom(E_0,E_0)_0\geqslant -\chi(E_0,E_0)_0=(r^2-1)(g-1).\]
The subscript $0$ means we are taking the traceless part.

Let a {\em Shatz stratum}
\begin{equation}\label{eq:shatz-str} A(r_1,d_1;\cdots ; r_{k},d_k)\subset A
\end{equation} be the locally closed subscheme of sheaves $E$, which have HN-type $(r_1,d_1;\cdots ; r_{k},d_k; r,d)$.
According to \cite{Sha77}, the HN-polygon rises when a sheaf specializes. Therefore, a higher HN-polygon corresponds to a smaller Shatz stratum.

We will mostly focus on {\em unstable Shatz strata}, namely strata parametrizing unstable sheaves of various HN-types. We will study the geometry of large unstable Shatz strata in the versal deformation space $A$ and the local topology of these strata along their boundaries.
We will mainly focus on higher rank cases.

For $r$, $d$, and $0<r_i<r$ fixed, we say $d_i$ is {\it minimal with respect to} $r_i$, or minimal, if
\begin{equation}
   \label{min-deg} d_i=\min\{n\in \mathbb{Z}\mid n/r_i>d/r\}.
\end{equation}
\begin{prop}\label{cod-unstable}
Suppose $r\geqslant 3$, $g\geqslant 3$, $r_1=s_i$ for some $1\leqslant i<t$, and $d_1$ is minimal with respect to $r_1$.
Then the closure $\overline{A(r_1,d_1)}$ contains $0$ and has pure codimension
\begin{equation*}
r_1(r-r_1)(g-1)-dr_1+d_1r
\end{equation*}
near $0$.
\end{prop}

When $E_0$ is a vector bundle, the result is nicer, see \cite[Corollary 15.4.3]{LeP97}. There, this minimality condition is not needed, and we also know the Shatz strata are smooth. This is because $\Ext^0_+$ is better controlled, \cite[Proposition 15.3.2]{LeP97}.

\begin{thm}\label{unstable-lci}
Suppose $r$, $g$, $r_1$, and $d_1$ are as in Proposition~\ref{cod-unstable}.
Let $E\in \ard$ be a vector bundle, which has Harder-Narasimhan filtration
\begin{equation}\label{hn}
0\subsetneqq F\subsetneqq E,
\end{equation}
such that $F$ and $E/F$ are semistable. Then, $\ard$ is a local complete intersection at $E$.
In particular, $\ard$ is locally irreducible at $E$.
\end{thm}

A special case of this result is when $E_0$ has HN-type $(r_1,d_1;r,d)$, then $A(r_1,d_1)$ contains $0$, which corresponds to $E_0$.
If we focus on the local irreducibility, we can relax the conditions on $g$ and $d_1$, see Proposition~\ref{smaller-strata}.

We next study the local topology of large unstable strata along their boundaries, when $E_0$ is not too far from being semistable.
\begin{thm}\label{loc-irr-1}
Suppose $r\geqslant 3$, $g\geqslant 2$, $0<r_1<r_2<r$, and $d_1$ and $d_2$ are minimal with respect to $r_1$ and $r_2$ respectively.
\begin{enumerate}[(a)]
    \item Assume $d_1/r_1=d_2/r_2$, $E_0$ has HN-type $(r_2,d_2;r,d)$, and there is only one filtration $0\subsetneqq F\subsetneqq E_0$ such that $\rk F=r_1$ and $\deg F=d_1$. If $\Hom(F,E_0/F)=0$, then $\overline{A(r_1,d_1)}$ is locally irreducible at $0$.
    \item Assume $(d-d_1)/(r-r_1)=(d-d_2)/(r-r_2)$, $E_0$ has HN-type $(r_1,d_1;r,d)$, and there is only one filtration $0\subsetneqq F\subsetneqq E_0$ such that $\rk F=r_2$ and $\deg F= d_2$. If $\Hom(F,E_0/F)=0$, then $\overline{A(r_2,d_2)}$ is locally irreducible at $0$.
\end{enumerate}
\end{thm}
See Figure~\ref{fig:two-strata}(I,II).
We note that the filtrations in the theorem are not HN-filtrations.
In either situation above, $E_0$ lies on the boundary of $\overline{A(r_1,d_1)}$ or $\overline{A(r_2,d_2)}$.

For an even more special $E_0$, we assume it has HN-type
\[(r_1,d_1;r_2,d_2;r,d)\]
with minimal $d_1$ and $d_2$.
The sheaf may appear in the intersection of closures of two maximal unstable strata, see \S\ref{intersection}.

\begin{thm}\label{loc-irr}
  Suppose $r\geqslant 3$, $g\geqslant 2$, and $d_1$ and $d_2$ are minimal.
  If $A_0\subset A$ is one of the Shatz strata $A(r_1,d_1)$, $A(r_2,d_2)$, or $A(r_1,d_1;r_2,d_2)$, then the closure $\overline{A_0}$ in $A$ is locally irreducible at $0$.
\end{thm}

The paper is organized as follows.
In \S\ref{deformations}, we study deformations of unstable sheaves to locally free stable sheaves and prove Proposition~\ref{cod-unstable} and Theorem~\ref{unstable-lci}. We also prove a Brill-Noether type result (Proposition \ref{brill-noether}), which is probably of general interests.
In \S\ref{sect:loc-irr}, we prove Theorems \ref{loc-irr-1} and \ref{loc-irr}.

\section{General deformations}\label{deformations}
We will prove Proposition~\ref{cod-unstable} and Theorem~\ref{unstable-lci} in thie section.

The next three lemmas concern general deformations of unstable vector bundles and sheaves with torsion.
\begin{lem}\label{stable-unstable}
Let $E$ be a coherent sheaf on $C$ with determinant $L$.
\begin{enumerate}[(a)]
\item \label{stable-unstable-1}
Suppose $E$ has torsion.
Then there is a deformation of $E$ whose general member is torsion free with determinant $L$.
\item \label{stable-unstable-2}
If $E$ is an unstable vector bundle, then there is a deformation of $E$ whose general member is stable with determinant $L$.
\item \label{torsion}
Suppose $E$ is $\mu$-unstable and has torsion.
Then there is a deformation of $E$ whose general member is $\mu$-unstable and locally free with determinant $L$.
\end{enumerate}
\end{lem}

\begin{proof}
(a) Let $E\cong F\oplus T$ where $F$ is the torsion free part and $T$ is the torsion part.
We can replace $T$ by $\oplus \oo_{x_i}$ where $x_i$'s are distinct points and $\oplus \oo_{x_i}$ has the same determinant as $T$.
We also have $\Ext^1(T,F)\cong \Hom(F,T\otimes \omega_C)^\vee\not=0$. Let $E^\prime\in \Ext^1(T,F)$ be a non-trivial extension.
If $E^\prime$ has torsion, let $T^\prime$ be the torsion.
There is an induced injective map $T^\prime \to T$.
It is not an isomorphism, otherwise it provides a splitting of the extension.
From the following commutative diagram,
\begin{equation*}
\begin{tikzcd}
  & & 0\arrow{d} & 0\arrow{d} &\\
&  & T^\prime\arrow{d} \arrow{r} & T^\prime \arrow{d}  &\\
0\arrow{r} & F\arrow{r}\arrow{d} & E^\prime \arrow{r}\arrow{d} & T\arrow{d} \arrow{r} & 0\\
0\arrow{r} & F\arrow{r} \arrow{d} & E^\prime/T^\prime \arrow{r} \arrow{d} & T/T^\prime  \arrow{r} \arrow{d} & 0\\
& 0 & 0 & 0 &
\end{tikzcd}
\end{equation*}
we know that all such extensions lie in the image of
$\Ext^1(T/T^\prime, F)\to \Ext^1(T,F)$.
The image is a proper subspace, according to the following exact sequence
\begin{equation*}
0\to \Ext^1(T/T^\prime, F)\to \Ext^1(T,F) \to \Ext^1(T^\prime,F)\to 0.
\end{equation*}
Since $T^\prime\subset T$ is a discrete datum, 
 there is a open subset of $\Ext^1(T,F)$ over which the restriction of the universal extension provides torsion free deformations of $E$.

(b) See \cite[Proposition 2.6]{NarRam75}.

(c) Let $G$ be the destabilizing quotient of $E$, which is locally free.
Then we have a short exact sequence
\begin{equation*}
0\to F\oplus T\to E\to G\to 0
\end{equation*}
where $T$ is torsion and $F$ is locally free.
Similar to the proof of (\ref{stable-unstable-1}), we can find a deformation $\{F_s\}_{s\in S}$ of $F\oplus T$ whose general member is locally free and unstable, and has the same determinant. Then, extensions of $G$ by $F_s$, trivial or not, provide the required deformation.
\end{proof}

We next show that an unstable vector bundle can be deformed to a ``minimally'' unstable bundle without traceless automorphisms, where the dimension of the deformation space can be calculated.

\begin{figure}[ht]
    \centering
    \includegraphics[scale=0.25]{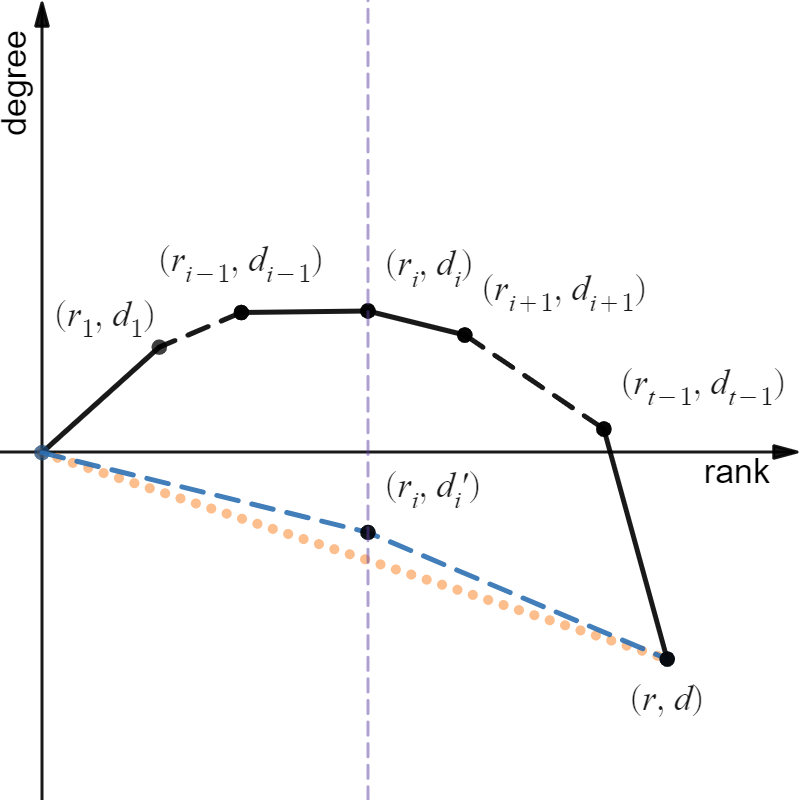}
    \caption{Deformation to a maximal unstable Shatz stratum}
    \label{fig:dec-deg}
\end{figure}

\begin{lem}\label{min-hn}
Suppose the genus $g\geqslant 3$.
Let $E$ be an unstable vector bundle on $C$ with $\rk(E)\geqslant 3$ and determinant $L$.
Suppose $E$ has Harder-Narasimhan type $(r_1,d_1;r_2,d_2;\ldots;r_{t-1},d_{t-1};r,d)$.
Then, for each pair $(r_i,d_i)$, $i=1,\ldots t-1$, there is a deformation of $E$ whose general member $E_s$ has determinant $L$ and no traceless automorphisms and fits into a short exact sequence of vector bundles
\begin{equation*}
0\to F_s\to E_s\to G_s\to 0.
\end{equation*}
where
\begin{enumerate}[(a)]
\item $F_s$ and $G_s$ are stable;
\item $F_s$ has rank $r_i$ and minimal degree $d_i^\prime$.
\end{enumerate}
\end{lem}
This lemma is illustrated by Figure~\ref{fig:dec-deg}.

\begin{proof}
The bundle $E$ fits in the following extension
\begin{equation*}
0\to F\to E\to G\to 0
\end{equation*}
such that $\rk F=r_i$ and $\deg F=d_i$.
We first show that we can decrease the degree of $F$. (In the process, the quotient will acquire torsion.) It is enough to decrease it by $1$ when $(\deg(F)-1)/\rk(F)>\mu(E)$.
Let $K$ be a sub-line bundle of $F$ with maximal degree \cite{Ati}. Then $F/K$ is a vector bundle, which can be assumed to be stable by Lemma~\ref{stable-unstable}.
Let $x$ be a point on $C$.
Then $F/(K(-x))\cong \oo_x\oplus F/K$.
Define a sub-vector bundle $F^\prime\subset F$ by the commutative diagram
\begin{equation*}
\begin{tikzcd}
0 \arrow{r} & K(-x) \arrow{r} \arrow[equal]{d} & F^\prime \arrow{r} \arrow[hook]{d} & F/K \arrow[hook]{d}\arrow{r} & 0\\
0\arrow{r} & K(-x) \arrow{r}  & F\arrow{r} & \oo_x\oplus F/K\arrow{r} & 0
\end{tikzcd}
\end{equation*}
where the second square is a pull-back square.
Then $F^\prime$ has degree $\deg(F)-1$.
Replace $F$ by $F^\prime$, which may be unstable.
We keep doing this until the degree of $F^\prime$ satisfies the following inequalities
\begin{equation*}
\rk(F)\mu(E)<\deg(F^\prime)\leqslant \rk(F)\mu(E)+1.
\end{equation*}
Namely, $\deg(F^\prime)$ is minimal.

By Lemma~\ref{stable-unstable}, we can assume $F^\prime $ to be stable and replace $E/F^\prime$ by a stable bundle, say $G^\prime$.
Furthermore, since $\rk (E)\geqslant 3$ and $g\geqslant 3$,
\begin{equation*}
\chi(G^\prime, F^\prime)=\rk(G^\prime )\rk(F^\prime)\left(\frac{\rk(E)}{\rk(G^\prime)\rk(F^\prime)}\left(\deg(F^\prime)-\rk(F^\prime)\mu(E)\right)+1-g\right)< 0,
\end{equation*}
we can choose a general $G^\prime$ such that
\begin{equation}\label{trivial-hom}
\Hom(G^\prime, F^\prime)=0,
\end{equation}
according to Proposition~\ref{brill-noether}.
Also, $\Ext^1(G^\prime, F^\prime)\not=0$.
We can replace $E$ by a non-trivial element $E^\prime$ in $\Ext^1(G^\prime, F^\prime)$.
We next show that $E^\prime$ has no traceless automorphism. We apply the functor $\Hom(-,F^\prime)$ to the short exact sequence
\begin{equation}\label{ses-generic-unstable}
0\to F^\prime\to E^\prime \to G^\prime\to 0.
\end{equation}
We have the associated long exact sequence
\begin{equation*}
0 \to \Hom(G^\prime, F^\prime)=0 \to \Hom(E^\prime, F^\prime) \to \Hom(F^\prime,F^\prime)\cong \mathbb{C}\to \cdots.
\end{equation*}
Since $E^\prime$ does not split, $\Hom(E^\prime, F^\prime)=0$.
Applying $\Hom(-,G^\prime)$ to (\ref{ses-generic-unstable}), we get an exact sequence
\begin{equation*}
0 \to \Hom(G^\prime, G^\prime)\cong \mathbb{C} \to \Hom(E^\prime, G^\prime) \to \Hom(F^\prime,G^\prime)=0.
\end{equation*}
Thus, $\Hom(E^\prime, G^\prime)\cong \mathbb{C}$. Therefore, from the long exact sequence
\begin{equation*}
0\to \Hom(E^\prime, F^\prime)=0\to \Hom(E^\prime, E^\prime) \to \Hom(E^\prime, G^\prime)\cong\mathbb{C}\to \cdots,
\end{equation*}
we conclude that $\Hom(E^\prime, E^\prime)_0=0$. This finishes the proof of the lemma.
\end{proof}

We have the following Brill-Noether type result.
\begin{prop}\label{brill-noether}
Let $C$ be a smooth projective algebraic curve of genus $g\geqslant 2$. Suppose $E$ and $F$ are general vector bundles over $C$ such that $\chi(E,F)\leqslant 0$. Then $\hom(E,F)=0$.
In particular, if $\chi(E,F)= 0$, then $\hom(E,F)=\ext^1(E,F)=
0$.
\end{prop}
\begin{proof}
The condition $\chi(E,F)\leqslant 0$ is equivalent to the condition on their slopes: $\mu(E)-\mu(F)\geqslant 1-g$. Since $E$ and $F$ are general and $g\geqslant 2$, we can assume them to be stable. 

Suppose there is a nontrivial morphism $\phi\colon E\to F$. Let $E^\prime=\ker \phi$ and $E^{\prime\prime}=\im \phi$. Let $r^\prime$, $d^\prime$, and $\mu^\prime$ denote the rank, degree, and slope of $E^\prime$. Symbols $r^{\prime\prime}$, $d^{\prime\prime}$, and $\mu^{\prime\prime}$ have similar meanings. Let $\mu_i^{\prime}$, $i=1,\dots, s$, denote the slopes of the Harder-Narasimhan gradings of $E^\prime$, and $\mu_j^{\prime\prime}$, $j=1,\dots, t$, denote those of $E^{\prime\prime}$. Then we have the inequalities
\begin{multline*}\mu_{\min}(E^\prime)=\mu_s^\prime<\dots<\mu_1^\prime=\mu_{\max}(E^\prime) <\mu(E)\\
<\mu_{\min}(E^{\prime\prime})=\mu_s^{\prime\prime}<\dots<\mu_1^{\prime\prime}=\mu_{\max}(E^{\prime\prime})<\mu(F).
\end{multline*}
Then $\hom(E^{\prime\prime},E^{\prime})=0$ and $\ext^1(E^{\prime\prime},E^{\prime})=-\chi(E^{\prime\prime},E^{\prime})=r^\prime r^{\prime\prime}(g-1)+r^\prime d^{\prime\prime}-r^{\prime\prime}d^\prime$.
The set of such $E$'s with nontrivial $\phi$ has dimension no greater than
\[(r^{\prime})^2(g-1)+1+(r^{\prime\prime})^2(g-1)+1+r^\prime r^{\prime\prime}(g-1)+r^\prime d^{\prime\prime}-r^{\prime\prime}d^\prime.\]
This number is strictly smaller than $\rk E^2(g-1)+1=\dim M^s(\rk E,\deg E)$, unless $\phi$ is injective. Here, $M^s(\rk E,\deg E)$ denotes the moduli space of stable vector bundles of $rk E$ and $\deg E$.

If $\phi$ is injective, a similar argument shows that the set of such $F$'s has dimension strictly smaller than $\dim M^s(\rk F,\deg F)$. We have thus proven the statement.
\end{proof}

We are now ready to prove Proposition~\ref{cod-unstable}.
\begin{proof}[Proof of Proposition~\ref{cod-unstable}]
According to Lemma \ref{min-hn} and Lemma~\ref{stable-unstable}(\ref{torsion}), the closure $\overline{\ard}$ contains $0$.
Take a general point $w\in \ard$ close to $0$.
To show the theorem, it suffices to the show that $\ard$ has the corresponding codimension at $w$.
Let $E_w$ denote the sheaf corresponding to $w$
and $(A^\prime,\mathscr{E}_{A^\prime};0^\prime, E_w)$ a versal deformation space of $E_w$.
Then there is an analytic neighborhood $U$ of $w\in A$ and an induced analytic map
\[f\colon (U, w)\to (A^\prime, 0^\prime).\]

Let $\urdo$ denote the set of isomorphism classes of unstable sheaves $E$ such that (i) $E$ has Harder-Narasimhan filtration of the form (\ref{hn}) where $F$ has rank $r_1$ and minimal degree $d_1$; (ii) $E$ has no traceless automorphisms. Notice that $d_1$ is determined by $r_1$, with $r$ and $d$ fixed.
We note that $f^{-1}(A^\prime \cap \urdo )=U\cap\urdo$, which follows from the definition of versal deformation space.
Here, the intersection $A^\prime \cap \urdo$ denotes elements in $A^\prime$ which also lie in $\urdo$. It is similar for the second intersection.
Hence, we only need to show:
\begin{enumerate}[(a)]
\item $f$ is submersive at $w$;
\item $\codim(A^\prime\cap \urdo, A^\prime)=r_1(r-r_1)(g-1)-dr_1+d_1r$.
\end{enumerate}
The argument on P.647 of \cite{Li97} is independent of the rank.
So, the same argument shows (a).
We next show (b).
Since $E_w$ has no traceless automorphism,
\begin{equation*}
\dim A^\prime=\ext^1(E_w,E_w)_0=(r^2-1)(g-1).
\end{equation*}
The dimension
\begin{align*}
\dim A^\prime\cap \urdo &=r_1^2(g-1)+1+(r-r_1)^2(g-1)+1-g+\ext^1(E_w/F, F)-1\\
&= (g-1)(r_1^2+r^2-rr_1-1)+dr_1-d_1r.
\end{align*}
Here, $F$ is as in (\ref{hn}). In the second equality, we have used the fact $\Hom(E_w/F, F)=0$. This is true for a general $E_w$, which has been shown in the proof of Lemma~\ref{min-hn}, (\ref{trivial-hom}).
A simple calculation gives the codimension.
\end{proof}

Equipped with Proposition~\ref{cod-unstable}, we can deduce Theorem~\ref{unstable-lci} quite easily.
\begin{proof}[Proof of Theorem~\ref{unstable-lci}]
As before, we use $\mathscr E_{A}$ to denote the versal family parametrized by $A$.
Let $\Quot_{\pi_A}(\mathscr E_{A},r-r_{1},d-d_{1})$ be the relative Quot scheme of quotients with rank $r-r_{1}$ and degree $d-d_{1}$. Here, $\pi_A$ is the projection $A\times C\to A$.
Let $a\in \ard$ denote the point corresponding to $E$.
Denote the quotient $E/F$ as $G$ and the quotient map as $q$.
Then we have the following exact sequence
\begin{equation*}
0\to \Hom(F,G)\to T_{q}\Quot(\mathscr E_{A}/A,r-r_{1},d-d_{1})\to T_{a}\ard \to \Ext^{1}(F,G).
\end{equation*}
By our assumption, $\Hom(F,G)=0$.
Thus, $\Ext^{1}(F,G)$ has dimension $-\chi(F,G)$, which is exactly the codimension of $\ard$ in $A$.
Therefore, $\ard$ is a local complete intersection at $a$.
With $r\geqslant 3$ and $g\geqslant 3$, $\dim A(r_1,d_1)>2$.
According to \cite[Theorem, P.185]{GorMac88} which is due to Hamm \cite{Ham71}, it is locally irreducible at $a$.
\end{proof}
\begin{rem}
The same argument actually can provide stronger local connectivity result.
\end{rem}

\section{Local (ir)reducibility of Shatz strata along boundaries}\label{sect:loc-irr}
In this section, we study versal deformation spaces of vector bundles, which are not too far from being semistable.
We will first prove Theorem \ref{loc-irr}, and then prove Theorem~\ref{loc-irr-1} in a very similar way.
Before these, we will review some basic facts about HN-polygons and list the possible relative positions of the lowest unstable HN-polygon involved.

Fix a coherent sheaf $E_0$ with rank $r$ and determinant $L$ on the curve $C$. Let $A$ be a versal deformation space of $E_0$ with fixed determinant.

\subsection{Harder-Narasimhan polygon under specialization}\label{intersection}
The HN-polygon rises when a sheaf specializes. But the converse is not necessarily true.
For example, we consider a complete family of vector bundles of rank 4 and degree 1 parametrized by a smooth scheme $S$.
The polygon $(1,2;4,1)$ lies above $(2,1;4,1)$, but the Shatz stratum corresponding to the first polygon has codimension $3g+4$ while the other one has codimension $4g-2$.
Moreover, the closure of a stratum is not necessarily a union of strata.

\begin{figure}[ht]
    \centering
    \begin{subfigure}[b]{0.45\textwidth}
    \includegraphics[width=\textwidth
    ]{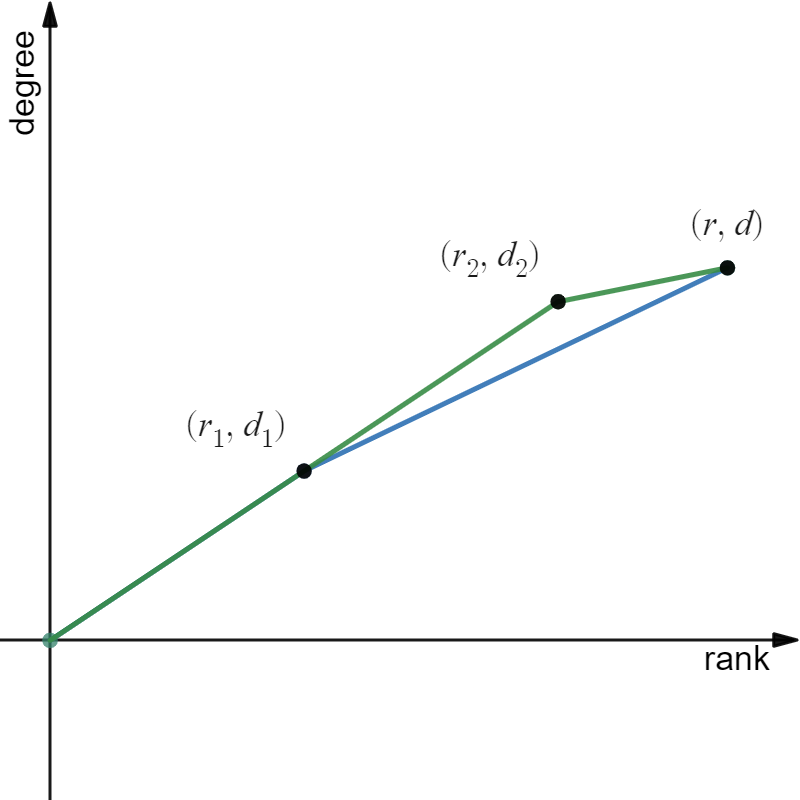}
    \caption{}
    \end{subfigure}\quad
    \begin{subfigure}[b]{0.45\textwidth}
    \includegraphics[width=\textwidth
    ]{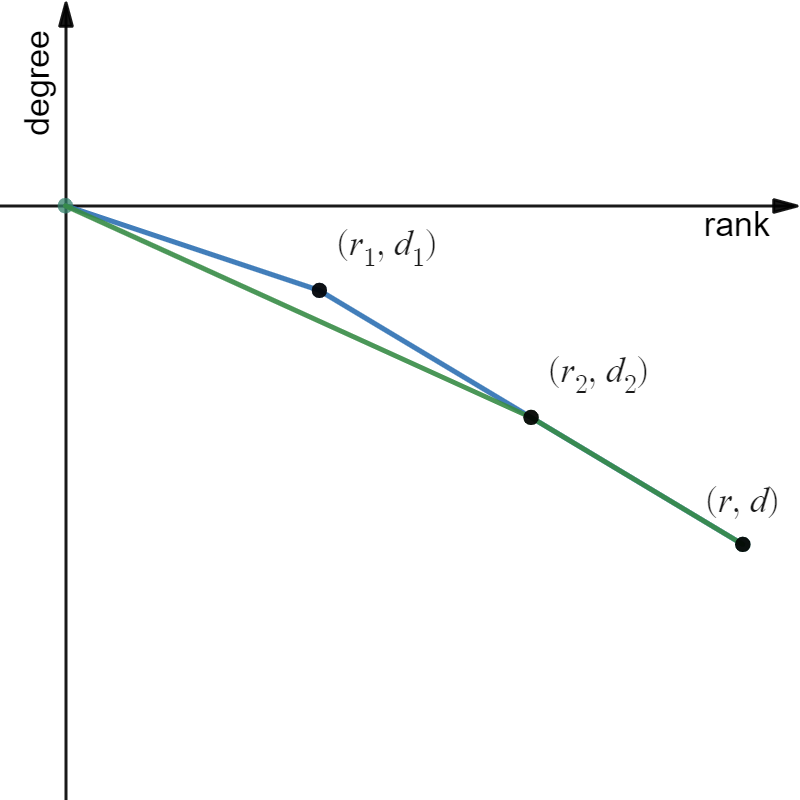}
    \caption{}
    \end{subfigure}\\
    \begin{subfigure}[b]{0.45\textwidth}
    \includegraphics[width=\textwidth
    ]{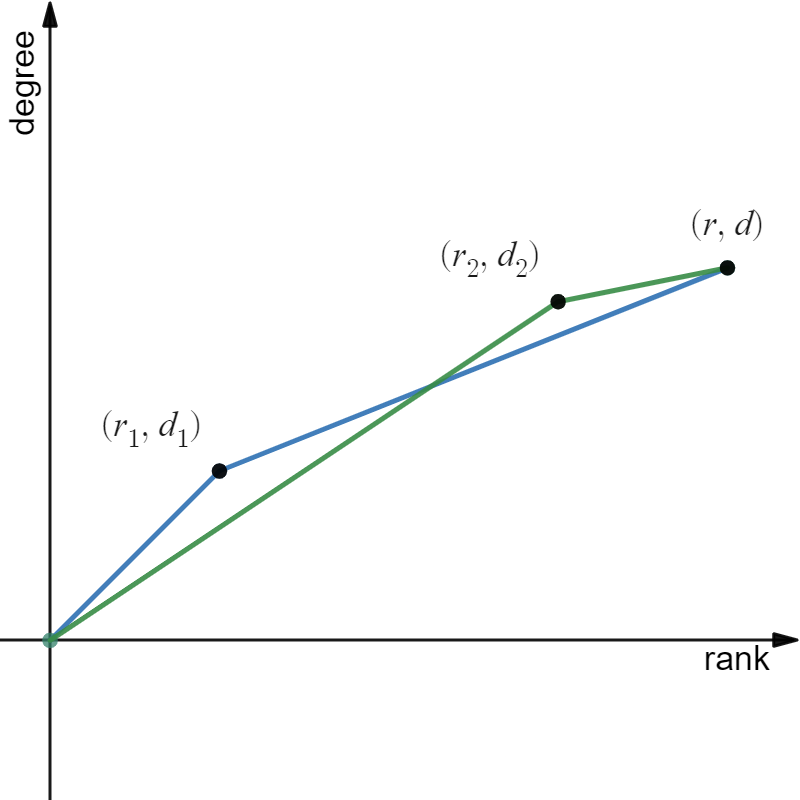}
    \caption{}
    \end{subfigure}\quad
    \begin{subfigure}[b]{0.45\textwidth}
    \includegraphics[width=\textwidth
    ]{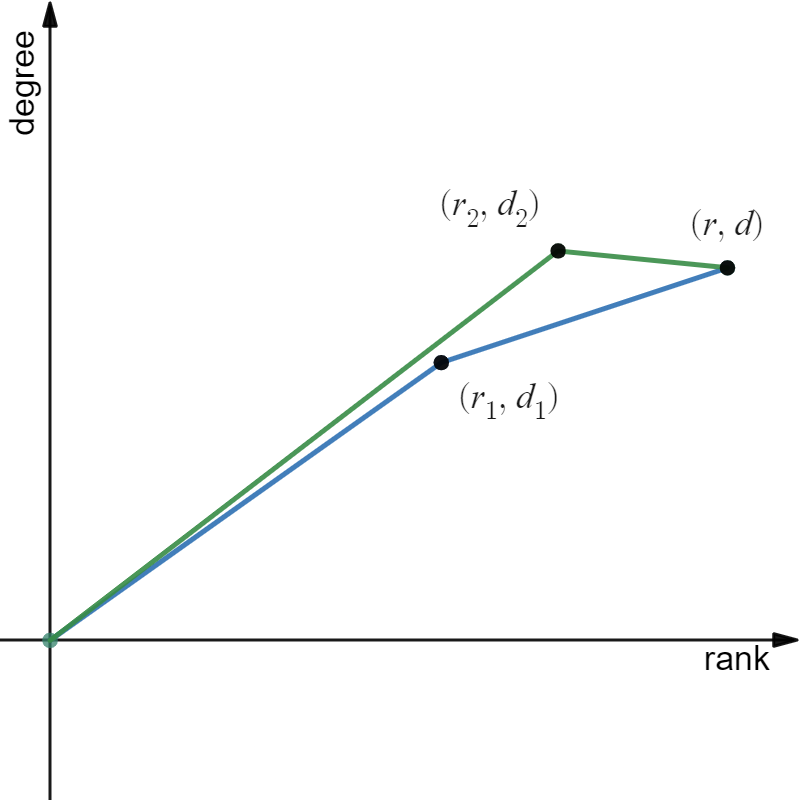}
    \caption{}
    \end{subfigure}
    \caption{Two maximal unstable Shatz strata}
    \label{fig:two-strata}
\end{figure}

Suppose we have $r_{1}$, $d_{1}$, $r_{2} $ and $d_{2}$ where $d_{1}$ and $d_{2}$ are minimal, namely, they satisfy the condition (\ref{min-deg}).
We define $\ard$ and $A(r_{2},d_{2})$ to be the locally closed subschemes of $A$ with the corresponding HN-types, as in (\ref{eq:shatz-str}).
Without loss of generality, suppose $r_{1}<r_{2}$.
By the minimality of $d_{1}$ and $d_{2}$, we are in one of the four cases in Figure~\ref{fig:two-strata}.

Numerically, each of the four cases can happen.
Except in Case (III), $(r_1,d_1;r_2,d_2;r,d)$ is not an HN-type.
We consider the intersection of $\overline{A(r_1,d_1)}$ and $\overline{A(r_2,d_2)}$ in these cases.
In Case (I), a general sheaf $E$ in the intersection fits in the following short exact sequence
\[0\to F\to E\to G\to 0\]
where $F$ is strictly semistable, $\rk F=r_2$, and $\deg F=d_2$.
In general, the boundary of $\overline{A(r_1,d_1)}$ contains a proper subset of $A(r_2,d_2)$. In particular, $\overline{A(r_1,d_1)}$ is not a union of strata.
In Case (II), a general sheaf $E$ in the intersection fits in a similar short exact sequence where $G$ is strictly semistable, $\rk G=r-r_1$, and $\deg G=d-d_1$.

\subsection{Proof of Theorem~\ref{loc-irr}}\label{subsect:loc-irr}

Assume that $E_0$ has HN-type
$(r_1,d_1;r_2,d_2;r,d)$ where $d_1$ and $d_2$ are minimal.
Then we are in Case (III) and $E_0$ lies in the intersection of $\overline{A(r_1,d_1)}$ and $\overline{A(r_2,d_2)}$.
Suppose the HN-filtration of $E_0$ is
\begin{equation}\label{3-step-fil}
    0\not=F_1\subsetneqq F_2\subsetneqq F_3=E_0.
\end{equation}

\subsubsection{
}\label{length-2}

Let $A_0=A(r_2,d_2)$ and $\overline{A_0}$ be its closure.
We consider the relative flag scheme $\Fl(\mathscr{E}, r_2,d_2)$ over $A$. The closed points are given by \[\{G_2\subsetneqq E\mid E\in A,\ \rk G_2=r_2, \ \deg G_2=d_2\}.\] Here, the flag scheme is actually a Quot scheme.
There is a forgetful map \[\pi\colon \Fl(\mathscr{E}, r_2,d_2)\to A,\] which is projective.
According to \cite[Propositions 15.4.1 \& 15.4.2]{LeP97}, (see also \cite{DreLeP85},) at the point $t=(F_2\subsetneqq E_0)$, we have an exact sequence
\begin{equation}
    0\to \Ext^0_{+}(E_0,E_0) \to T_{t}\Fl(\mathscr{E}, r_2,d_2) \xrightarrow{T_t\pi} T_0A\cong \Ext^1(E_0,E_0)_0\xrightarrow{\omega_+} \Ext^1_{+}(E_0,E_0)\to 0.
\end{equation}
The $\Ext^i_+$ here is calculated with respect to the filtration $F_2\subsetneqq E_0$.
Notice that $\overline{A_0}$ is the image of $\pi$ and \begin{equation}\label{only-fil}
    \mbox{there is only one point }t \mbox{ lying above }0\in A.
\end{equation}
There is a spectral sequence converging to $\Ext^{p+q}_+(E_0,E_0)$:
\begin{equation}\label{ext+ss}
    E_1^{p,q}=\begin{cases}
    \bigoplus_i\Ext^{p+q}(\gr_i(E_0), \gr_{i-p}(E_0) ) & \text{if $p<0$,} \\
   0 & \text{otherwise.}
    \end{cases}
\end{equation}
In this case,
\begin{equation}\label{eq:no-hom-fil}
E^{-1, 1}_1=\Hom(\gr_1(E_0),\gr_2 (E_0))=\Hom (F_2,E_0/F_2)=0.\end{equation}
The last equality is obtained via semistability and slope comparison. Thus, $\Ext^0_+(E_0, E_0)=0$.
On the othe hand,  $E^{-1,2}_1=\Ext^1(F_2,E_0/F_2)$, and other terms are zero, since we are working with a 2-step filtration over a curve.
Therefore, the spectral sequence degenerates and
\begin{equation*}
    \ext^1_+(E_0,E_0)=\ext^1(F_2,E_0/F_2)=-\chi(F_2,E_0/F_2)=d_2r-r_2d+r_2(r-r_2)(g-1).
\end{equation*}
According to Mori \cite[\S1]{Mor79} and Li \cite{Li94-kodaira}, the number of defining equations of $\overline{A_0}\subset A$ at $0$ is bounded by dimension of the obstruction space $\Ext^1_{+}(E_0,E_0)$, which is the number above. (See also \cite[Proposition 2.A.11]{HL10}.)
Applying Hamm's result again, $\overline{A_0}$ is locally irreducible, as long as \begin{equation*}
    \dim A-\ext^1_+(E_0,E_0)-2\geqslant 0.\end{equation*}
This is satisfied when $g\geqslant 2$ and $r\geqslant 3$. Here, we have used the assumption $(d_2-1)/r_2\leqslant d/r$.

When $A_0=A(r_1,d_1)$, 
we can prove the local irreducibility of $\overline{A_0}$ at $0$ in the same way.

\subsubsection{
}
When $A_0=A(r_1,d_1;r_2,d_2)$, we consider the relative flag scheme $\Fl(\mathscr{E};r_1,d_1;r_2,d_2)$ of filtrations of the form $G_1\subsetneqq G_2\subsetneqq E$ with the corresponding ranks and degrees. Then, similarly, the number of defining equation of $\overline{A_0}\subset A$ at $0$ is still given by $\ext^1_+(E_0,E_0)$ with respect to the filtration (\ref{3-step-fil}).
In the spectral sequence (\ref{ext+ss}),
\begin{align*}
    E^{-1,1}_1&=\Hom(\gr_1 E_0,\gr_2 E_0)\oplus \Hom(\gr_2 E_0,\gr_3 E_0),\\
    E^{-1,2}_1& =\Ext^1(\gr_1 E_0,\gr_2 E_0)\oplus \Ext^1(\gr_2 E_0,\gr_3 E_0),\\
    E^{-2,2}_1&=\Hom(\gr_1 E_0,\gr_3 E_0),\  \mbox{and}\\
    E^{-2,3}_1&=\Ext^1(\gr_1 E_0,\gr_3 E_0),
\end{align*}
while other terms are zero.
Moreover, $E^{-1,1}_1=E^{-2,2}_1=0$ due to semistability and slope comparison. Therefore, the spectral sequence degenerates as well in this case, and
\begin{eqnarray*}
    \ext^1_+(E_0,E_0)&=&\sum_{1\leqslant i<j\leqslant 3}\ext^1(\gr_i E_0,\gr_j E_0)=-\sum_{1\leqslant i<j\leqslant 3}\chi(\gr_i E_0,\gr_j E_0)\\
    &=& r_2d_1-r_1d_2+r_3d_2-r_2d_3+(r_1r_2-r_1^2+r_2r_3-r_2)(g-1).
\end{eqnarray*}
When $g\geqslant 2$ and $r\geqslant 3$, we also have the local irreducibility of $\overline{A_0}$ at $0$.

\subsection{Proof of Theorem~\ref{loc-irr-1}}
\begin{proof}
\begin{enumerate}[(a)]
    \item This corresponds to Case (I) in Figure~\ref{fig:two-strata}.
    Notice that $F$ is a semistable and $E_0/F$ has HN-type $(r_2-r_1,d_2-d_1;r-r_1,d-d_1)$.
    We replace the filtration (\ref{3-step-fil}) by the following filtration
\begin{equation}\label{eq:non-hn-fil}
    0\not=F_1=F\subsetneqq F_2\subsetneqq F_3=E_0,
\end{equation}
such that $F_2/F\subsetneqq E_0/F$ is the HN-filtration of $E_0/F$.
Note that (\ref{eq:non-hn-fil}) is not the HN-filtration.
Then the assumptions imply that conditions (\ref{only-fil}) and (\ref{eq:no-hom-fil}) still hold and the argument in \S\ref{subsect:loc-irr} applies here as well.
\item This corresponds to Case (II) in Figure~\ref{fig:two-strata}. The argument is similar. We only need to replace (\ref{3-step-fil}) by
    \begin{equation*}
    0\not=F_1\subsetneqq F_2=F\subsetneqq F_3=E_0,
\end{equation*}
such that $F_1\subsetneqq F$ is the HN-filtration of $F$.
\end{enumerate}

\end{proof}
\subsubsection{Locally reducible examples}
Counterexamples appear when the condition (\ref{only-fil}) does not hold.

Assume $E_0$ has HN-type $(1, d_1;3,d)$ such that $d-d_1$ is even. Furthermore, we assume it fits in the following short exact sequence
\begin{equation*}
    0\to F\to E_0 \to L_1\oplus L_2 \to 0
\end{equation*}
such that $L_1$ and $L_2$ are line bundles of degree $(d-d_1)/2$, which are not isomorphic.
We will explain that the subset $\overline{A(2, (d+d_1)/2)}\subset A$ is locally reducible at $0$.

Let $d_2= (d+d_1)/2$. Let $\mathscr{E}\to A\times C$ be the associated versal family.
We consider the relative flag scheme \[\Fl=\Fl(2,d_2;3,d) \to A\] parametrizing filtrations of the type $(2,d_2;3, d)$.
Then, $\Fl$ is mapped to $\overline{A(2, d_2)}$.
For $i=1$ or $2$, we define and consider the short exact sequence
\begin{equation*}
    0\to K_i\to E_0 \to L_i\to 0.
\end{equation*}
We denote the corresponding filtration as $s_i\in \Fl$.
 Then it gives rise to an exact sequence
 \begin{equation*}
     0\to \Ext^0_{+,i}(E_0,E_0) \to T_{s_i}\Fl \to T_0A\cong \Ext^1(E_0,E_0)_0\to \Ext^1_{+,i}(E_0,E_0).
 \end{equation*}
There is a converging spectral sequence 
$E_1^{p,q}\Rightarrow \Ext^{p+q}_{+,i}(E_0,E_0)$
with $E_1^{-1,1}= E_1^{-2,2} =E_1^{-2,3}=0$ and $E_1^{-1,2}=\Ext^1(K_i,L_i)$.
Then we have an exact sequence
\[0\to T_{s_i}\Fl \to \Ext^1(E_0,E_0) \to \Ext^1(K_i,L_i)\to 0.\]
It is indeed exact on the right according to \cite[Proposition 15.4.1]{LeP97}.
The map $\Ext^1(E_0,E_0) \to \Ext^1(K_i,L_i)$ is the natural one induced by the inclusion $K_i\hookrightarrow E_0 $ and the surjective map $E_0\twoheadrightarrow L_i$, namely, it is the composition:
\[\Ext^1(E_0,E_0)\to \Ext^1(E_0, L_i)\to \Ext^1(K_i,L_i).\]

Notice we have the short exact sequence $0\to F\to K_1\to L_2\to 0$.
Since $\Hom(F,L_1)=0$, we have an inclusion $\Ext^1(L_2,L_1)\cong \mathbb{C}^g\hookrightarrow \Ext^1(K_1, L_1)$.
We have the following commutative diagram.
\begin{equation*}
    \begin{tikzcd}
    \Ext^1(L_2, K_2)\arrow[two heads]{rr} \arrow{d} & & \Ext^1(L_2, L_1)\arrow[hook]{d}\\
    \Ext^1(E_0, K_2)\arrow{r}{\phi} & \Ext^1(E_0, E_0)\arrow{r} & \Ext^1(K_1, L_1).
    \end{tikzcd}
\end{equation*}
Here, all the maps are the obvious ones.
The composition of the bottom row is nontrivial.
On the other hand, $\phi$ factors through $T_{s_2}\Fl$.
Therefore, $T_{s_2}\Fl \to \Ext^1_{+,1}(E_0,E_0)$ is nontrivial.
We can prove in the same way $T_{s_1}\Fl \to \Ext^1_{+,2}(E_0,E_0)$ is nontrivial.
Hence, $\overline{A(2, d_2)}\subset A$ is locally reducible at $0$.

This example does not contradict Theorem~\ref{unstable-lci}. The theorem is about the stratum, while this example demonstrates that local reducibility may develop along the boundary of a stratum. In the theorem, $A(r_1,d_1)$ may contain $0$.

\subsection{Smaller Shatz strata} An argument similar to \S\ref{length-2} can show the following statement.
\begin{prop}\label{smaller-strata}
Suppose $g\geqslant 2$, $E_0$ has HN-type $(r_1,d_1;r,d)$, and
\[(r^2-1-r_1r+r_1^2)(g-1)-d_1r+dr_1\geqslant 2.\]
Then $A(r_1,d_1)$ is locally irreducible at $0$.
\end{prop}

\vskip6pt
{\em Acknowledgment.} The author would like to thank Jun Li, Howard Nuer, and Xiaolei Zhao for helpful discussions. The author also would like to express his gratitude towards Shanghai Center for Mathematical Sciences and University of California, Santa Babara for their hospitality during his visits.

\bibliographystyle{alpha}
\bibliography{picard}

\begin{thebibliography}{Ham71}

\bibitem[Ati57]{Ati}
M.~Atiyah.
\newblock Vector bundles over an elliptic curve.
\newblock {\em Proc. London Math. Soc. (3)}, 7:414--452, 1957.

\bibitem[DLP85]{DreLeP85}
J.-M. Drezet and J.~Le~Potier.
\newblock Fibr\'{e}s stables et fibr\'{e}s exceptionnels sur {${\bf P}_2$}.
\newblock {\em Ann. Sci. \'{E}cole Norm. Sup. (4)}, 18(2):193--243, 1985.

\bibitem[GM88]{GorMac88}
Mark Goresky and Robert MacPherson.
\newblock {\em Stratified {M}orse theory}, volume~14 of {\em Ergebnisse der
  Mathematik und ihrer Grenzgebiete (3) [Results in Mathematics and Related
  Areas (3)]}.
\newblock Springer-Verlag, Berlin, 1988.

\bibitem[Ham71]{Ham71}
Helmut Hamm.
\newblock Lokale topologische {E}igenschaften komplexer {R}\"{a}ume.
\newblock {\em Math. Ann.}, 191:235--252, 1971.

\bibitem[HL10]{HL10}
D.~Huybrechts and M.~Lehn.
\newblock {\em The geometry of moduli spaces of sheaves}.
\newblock Cambridge Mathematical Library. Cambridge University Press,
  Cambridge, second edition, 2010.

\bibitem[Li94]{Li94-kodaira}
Jun Li.
\newblock Kodaira dimension of moduli space of vector bundles on surfaces.
\newblock {\em Invent. Math.}, 115(1):1--40, 1994.

\bibitem[Li97]{Li97}
Jun Li.
\newblock The first two {B}etti numbers of the moduli spaces of vector bundles
  on surfaces.
\newblock {\em Comm. Anal. Geom.}, 5(4):625--684, 1997.

\bibitem[LP97]{LeP97}
J.~Le~Potier.
\newblock {\em Lectures on vector bundles}, volume~54 of {\em Cambridge Studies
  in Advanced Mathematics}.
\newblock Cambridge University Press, Cambridge, 1997.
\newblock Translated by A. Maciocia.

\bibitem[Mor79]{Mor79}
Shigefumi Mori.
\newblock Projective manifolds with ample tangent bundles.
\newblock {\em Ann. of Math. (2)}, 110(3):593--606, 1979.

\bibitem[NR75]{NarRam75}
M.~S. Narasimhan and S.~Ramanan.
\newblock Deformations of the moduli space of vector bundles over an algebraic
  curve.
\newblock {\em Ann. Math. (2)}, 101:391--417, 1975.

\bibitem[Sha77]{Sha77}
Stephen~S. Shatz.
\newblock The decomposition and specialization of algebraic families of vector
  bundles.
\newblock {\em Compositio Math.}, 35(2):163--187, 1977.

\end{thebibliography}
\end{document}